\documentclass[]{aspm}

\articleinfo{}{}{}

\usepackage{verbatim}
\usepackage{amssymb}
\usepackage{amsbsy}
\usepackage{amscd}
\usepackage{amsmath}
\usepackage{amsthm}
\usepackage[mathscr]{eucal}
\usepackage{amsmath,amssymb,latexsym, geometry, amsthm, amsfonts, multicol, color, graphicx}

\newtheorem{theorem}{Theorem}
\newtheorem{defn}{Definition}\numberwithin{defn}{section}
\newtheorem{example}[theorem]{Example}
\newtheorem{proposition}[theorem]{Proposition}
\newtheorem{fact}[theorem]{Fact}
\newtheorem{remark}[theorem]{Remark}
\newtheorem{question}[theorem]{Question}

\title[Billey's formula]{Billey's formula in combinatorics, geometry, and topology}

\author[J.S. Tymoczko]{Julianna Tymoczko}

\address{Department of Mathematics, Smith College, 44 College Lane, Northampton MA 01063, USA}

\email{jtymoczko@smith.edu}

\rcvdate{}
\rvsdate{}



\thanks{The author wishes to thank the organizers of the International Summer School and Workshop on Schubert Calculus for the invitation to present this paper, and the NSF for the support of grants DMS-1248171 and DMS-1205283.}

\begin{document}

\maketitle

\section{Introduction}

In this paper we describe a powerful combinatorial formula and its implications in geometry, topology, and algebra.  This formula first appeared in the appendix of a book by Andersen, Jantzen, and Soergel \cite[Appendix D]{AJS94}.  Sara Billey discovered it independently five years later, and it played a prominent role in her work to evaluate certain polynomials closely related to Schubert polynomials \cite{Bil99}.

To set the stage for our discussion, we review well-known foundations of Schubert calculus in Lie type $A_{n-1}$.  Consider the group of invertible matrices $GL_n(\mathbb{C})$ with the subgroup $B$ of upper-triangular matrices.  The flag variety is the quotient $GL_n(\mathbb{C})/B$ and can  be thought of geometrically as the collection of nested vector subspaces $V_1 \subseteq V_2 \subseteq \cdots V_{n-1} \subseteq \mathbb{C}^n$ where each $V_i$ is $i$-dimensional.  The geometry of the flag variety is interwoven with the combinatorics of the permutation group: the torus $T$ of diagonal matrices in $B$ acts on the flag variety, and its fixed points are the flags corresponding to permutation matrices.  For each permutation $w$, the double coset $BwB$ is an affine cell inside the flag variety, and the union of these double cosets forms a CW-decomposition.  The closures of the cells $\overline{BwB}$ are the Schubert varieties, which induce a basis for the cohomology of the flag variety.  Combinatorial properties of the permutations $w$ determine topological properties of Schubert varieties: for instance, the number of inversions of $w$ counts the dimension of the variety $\overline{BwB}$ in the flag variety.   

Billey's formula relates all of these pieces of Schubert calculus: the geometry of Schubert varieties, the action of the torus $T$ on the flag variety, combinatorial data about permutations, the cohomology of the flag variety and of the Schubert varieties, and the combinatorics of root systems (the generalization of inversions of a permutation).

Combinatorially, Billey's formula describes an invariant of pairs of elements of a Weyl group.  On its face, this formula is a combination of roots built from subwords of a fixed word.  As we will see, it has deeper geometric and topological meaning as well:
\begin{itemize}
\item It tells us about the tangent spaces at each permutation flag in each Schubert variety.
\item It tells us about singular points in Schubert varieties.
\item It tells us about the values of Kostant polynomials.
\end{itemize}
Billey's formula also reflects an aspect of GKM theory, which is a way of describing the torus-equivariant cohomology of a variety just from information about the torus-fixed points in the variety.

We will also describe some applications of Billey's formula, including concrete combinatorial descriptions of Billey's formula in special cases, and ways to bootstrap Billey's formula to describe the equivariant cohomology of subvarieties of the flag variety to which GKM theory does not apply.

More precisely, Section \ref{section: combinatorics} states Billey's formula, gives examples, and lists its main properties.  Section \ref{section: motivation} describes the motivation behind Billey's formula and the various places it arose.  Sections \ref{section: applications} and \ref{section: poset pinball} then describe newer applications of Billey's formula: a tool called {\em excited Young diagrams} with which one can compute Billey's formula for Grassmannians, and a tool called {\em poset pinball} with which one can use Billey's formula to compute the equivariant cohomology of various subvarieties of the flag variety.  Section \ref{section: Schubert calculus for Petersons} even shows how poset pinball can be used to construct a complete Schubert calculus for the Peterson variety, a singular subvariety of the flag variety with important applications to quantum cohomology.  Finally Section \ref{section: questions} concludes with some open questions and conjectures.

We used the example of Lie type $A$ in this introduction, and the reader is free to consider that example throughout.  However, the theorems in this paper apply to all complex linear algebraic groups.  We use the following notation throughout the manuscript.
\begin{itemize}
\item $G$ is a complex linear algebraic group
\item $B$ is a Borel subgroup of $G$ 
\item $T$ is a maximal torus in $B$
\item $\Phi$ is the root system corresponding to $G$ 
\item $\Delta=\{\alpha_1, \ldots, \alpha_n\}$ is the set of positive simple roots corresponding to $B$
\item $W$ is the Weyl group associated to $G$ and $T$
\item All cohomology is taken with complex coefficients
\end{itemize}

\section{Billey's Formula in Combinatorics}\label{section: combinatorics}

We begin by introducing Billey's formula from a strictly combinatorial viewpoint.  From this point of view, Billey's formula  simply associates a polynomial to each pair of elements in the Weyl group.  More formally, for each pair $v, w$ in the Weyl group $W$,
Billey's formula describes a polynomial $\sigma_v(w)$ 
in the simple roots $\alpha_1, \alpha_2, \ldots, \alpha_n$. 

Billey's formula uses two main tools:
\begin{itemize}
\item A {\em reduced word} $b_1b_2b_3 \cdots b_m$ is an expression for which each $b_i$ is a simple reflection in $W$ and the product $b_1b_2b_3 \cdots b_m$ cannot be written in fewer simple reflections.
\item To each reduced word $b_1b_2b_3 \cdots b_m$ and integer $j$ such that $1 \leq j \leq m$ we associate the root $r_j= b_1b_2b_3 \cdots b_{j-1}(\alpha_j)$.
\end{itemize}
Reduced words are commonly used when considering Weyl groups, Coxeter groups, or indeed any group with prescribed generator sets.  

The roots $r_j$ deserve more comment: they are important in Lie theory, though mysterious to some topologists.  The roots $r_j$ are always positive \cite[Section 10.2]{Hum78}.  They are the roots negated by $w$, in the sense that every $r_j$ is $-w(\alpha)$ for some positive root $\alpha$.  Even more important, the set of roots $r_j$ determine $w$.  More precisely, let $b_1b_2b_3 \cdots b_m$ be a reduced word for $w$ and denote by $N(w)$ the set $\{r_1, r_2, \ldots, r_m\}$ of all positive roots that can be formed in this way.  It turns out that the set $N(w)$ uniquely specifies $w$ in the sense that if $N(w_1) = N(w_2)$ then $w_1 = w_2$ \cite{Kos61}.  These sets $N(w)$ satisfy many other interesting combinatorial properties.  %Moreover, these sets $N(w)$ satisfy other useful combinatorial properties, like being closed under addition (in the sense that if $\alpha, \beta \in N(w)$ and $\alpha+\beta$ is a root then $\alpha+\beta \in N(w)$) \cite{???}.  The interested reader is encouraged to learn more about these sets \cite{???}; the casual reader can simply rest assured that the roots $r_j$ are natural from an algebraic and combinatorial viewpoint.

\begin{defn}[Billey's formula \cite{Bil99, AJS94}] 
Let $v,w$ be elements of $W$.  Fix a reduced word $b_1b_2b_3\cdots b_m$ for $w$.  Billey's formula for $\sigma_v(w)$ is the polynomial
\[\sigma_v(w) = \sum  \prod_{i=1}^k  r_{j_1} r_{j_2} \cdots r_{j_k}\]
where the sum is over all reduced subwords $b_{j_1}b_{j_2} \cdots b_{j_k}$ for $v$ in $w$, and for each $j$
\[r_j = b_1b_2b_3 \cdots b_{j-1}(\alpha_j).\]
\end{defn}
Andersen, Jantzen, and Soergel gave this first, in a remark to be proven by the reader  \cite[Appendix D, Remark on p. 298]{AJS94}.

\begin{example}\label{example: first Billey example}
Let $v=s_1$ and $w=s_2s_1s_2$ in the symmetric group $S_3$, namely the Weyl group of type $A_2$. 

\begin{enumerate}
\item[\em Step 1:] Find all possible subwords of $w$ that equal $v$:
\begin{center}
\begin{picture}(15,5)(30,0)
\put(0,0){$s_2$}
\put(15,0){$s_1$}
\put(19,1){\circle{12}}
\put(30,0){$s_2$}
\end{picture}
\end{center}
\item[\em Step 2:] Compute the roots in each term:
\begin{center}
\begin{picture}(25,5)(50,0)
\put(0,0){$s_2$}
\put(15,0){$\alpha_1$}
\put(19,1){\circle{13}}
\put(30,0){\color{green} $s_2$}
\put(40,0){$=\alpha_1+\alpha_2$}
\end{picture}
\end{center}
\item[\em Step 3:] Add terms:
\[\sigma_{s_1}(s_2s_1s_2) = \alpha_1+\alpha_2.\]
\end{enumerate}
Billey's formula in this case is the polynomial $\alpha_1+\alpha_2$.
\end{example}

We now consider several extreme cases of the definition.

First, suppose $v$ is not a subword of $w$.  In this case, the first step of the algorithm (``find all possible subwords...")  fails and there are no terms in Billey's formula.  More formally, we have the following.
\begin{proposition}\label{proposition: upper-triangular}
If $v \not \leq w$ in Bruhat order then $\sigma_v(w)=0$.
\end{proposition}

Second, suppose $v$ is the identity $e \in W$.  In this case, every element $w \in W$ contains $e$ trivially and uniquely as the empty subword.  We take the root computed in the second step of the algorithm to be $1$, giving us the following.
\begin{proposition}
If $v=e$ then $\sigma_e(w)=1$ for all $w \in W$p.
\end{proposition}

Finally, suppose $v$ and $w$ are both the longest element $w_0 \in W$.  In this case, Billey's formula gives us a well-known quantity: the product of the positive roots.  (See also the subset $N(w_0)$ above.)
\begin{proposition}
If $v=w=w_0$ is the longest element in the Weyl group then
\[\textup{FACT:          } \hspace{0.15in} \sigma_{w_0}(w_0) = \prod_{\alpha \in \Phi^+} \alpha.\]
\end{proposition}

Other important properties of Billey's formula are more subtle; before discussing them, we give one more substantive example.

\begin{example}
Let $v=s_1$ and $w=s_1s_2s_1$ in type $A$. 

\begin{enumerate}
\item[\em Step 1:] Find all possible subwords of $w$ that equal $v$:
\begin{center}
\begin{picture}(15,5)(30,0)
\put(0,0){$s_1$}
\put(4,1){\circle{12}}
\put(15,0){$s_2$}
\put(30,0){$s_1$}
\end{picture}
 \hspace{0.25in} and \hspace{0.25in}
\begin{picture}(15,5)(0,0)
\put(0,0){$s_1$}
\put(15,0){$s_2$}
\put(30,0){$s_1$}
\put(34,1){\circle{12}}
\end{picture}
\end{center}
\item[\em Step 2:] Compute the roots in each term:
\begin{center}
\begin{picture}(15,5)(30,0)
\put(0,0){$\alpha_1$}
\put(5,1){\circle{14}}
\put(15,0){\color{green} $s_2$}
\put(30,0){\color{green} $s_1$}
\end{picture}
 \hspace{0.25in} and \hspace{0.25in}
\begin{picture}(15,5)(0,0)
\put(0,0){$s_1$}
\put(15,0){$s_2$}
\put(30,0){$\alpha_1$}
\put(35,1){\circle{14}}
\end{picture}
\end{center}
\item[\em Step 3:] Add terms:
\[\sigma_{s_1}(s_1s_2s_1) = \alpha_1+\alpha_2.\]
\end{enumerate}
Using a different word for $w$ than Example \ref{example: first Billey example}, Billey's formula still gives $\alpha_1+\alpha_2$.
\end{example}

This leads directly to the fundamental properties of $\sigma_v(w)$.  Recall that the length $\ell(w)$ of a Weyl group element $w \in W$ is the number of simple reflections in a reduced word for $w$, namely if $w = b_1 b_2 \cdots b_m$ then $\ell(w)=m$.

\begin{proposition}
Let $v,w \in W$.  Billey's formula $\sigma_v(w)$ is:
\begin{itemize}
\item a polynomial in $\alpha_1, \alpha_2, \ldots, \alpha_n$ with nonnegative integer coefficients
\item a polynomial of degree $\ell(v)$
\item independent of the choice of reduced word for $w$
\end{itemize}
\end{proposition}
The first claim in this proposition follows from the fact that each root $r_j$ is always positive, and that every positive root can be written in terms of the simple roots with nonnegative integer coefficients.  The second claim follows from the algorithm itself: each term in $\sigma_v(w)$ is the product of roots, one for each simple reflection in a reduced word for $v$.  Hence each term has degree $\ell(v)$.  The third claim is not trivial; Billey proved it in her original manuscript \cite{Bil99}, and Anderson, Jantzen, and Soergel assert the claim without proof \cite{AJS94}.

\section{What Billey's formula means} \label{section: motivation}

We now consider the motivation behind the discovery of Billey's formula and other places it first emerged.

\subsection{Representations of quantum groups}

Andersen, Jantzen, and Soergel wanted to answer a problem of a different flavor from the Schubert calculus described in this manuscript: identifying representations of quantum groups and of semisimple groups over fields of characteristic $p$.  Their approach was to fit their specific questions into a larger framework.  They built a more abstract endomorphism algebra and showed that its properties captured their original representation-theoretic questions.

As a small application of their results, Andersen, Jantzen, and Soergel consider a particular endomorphism ring studied earlier by Kostant and Kumar \cite{KosKum86, KosKum90}.  They proved that classes $\sigma_w$ constructed by Kostant and Kumar actually generate various sub- and quotient modules of the endomorphism ring. In order to do this, Andersen, Jantzen, and Soergel explicitly identified the polynomials $\sigma_v(w)$ \cite[Chapter 19, Appendix D]{AJS94}.

\subsection{Orbit values of Kostant polynomials}
Billey's original goal was to study the values of certain polynomials called {\em Kostant polynomial} that are also related to the classes $\sigma_w$.  Let $O$ denote a regular element of the torus $\mathfrak{t}$.  Kostant polynomials are (nonhomogeneous) elements of $\mathbb{C}[\mathfrak{t}^*]$ parametrized by the Weyl group, of degree determined by the length of the Weyl group element associated to the polynomial, and defined by certain vanishing conditions on the orbits $Ov$ for each $v \in W$.   Surprisingly, these polynomials are essentially unique: Kostant showed that they are unique up to the ideal in $\mathbb{C}[\mathfrak{t}^*]$ that vanishes on the orbits $OW$.  Moreover, the highest homogeneous component is a Schubert class in $H^*(G/B)$. 

Later work expanded this deep connection between the Kostant polynomials and the cohomology of the flag variety. Carrell proved that the cohomology $H^*(G/B)$ is isomorphic as a graded ring to the coordinate ring of the variety associated to the points $OW$ \cite{Car89}.  Kostant and Kumar generalized Kostant polynomials to the collections of polynomials $\sigma_w$ and in a series of papers proved that the classes $\sigma_w$ generate the equivariant cohomology of the flag variety $G/B$ \cite{KosKum86, KosKum90}.

\subsection{Kumar's criterion for Schubert varieties}

Of course, algebraic information about Schubert classes is closely related to geometric data of Schubert varieties.  Kumar pursued this relationship, strengthening his analysis of some specific polynomials $\sigma_v(w)$ to give a criterion for when the (opposite) Schubert variety $X^v$ is smooth at the flag $wB$.  It turns out that smoothness is equivalent to $\sigma_v(w)$ being a specific product of distinct positive roots \cite{Kum96}.

\begin{theorem}[Kumar's Criterion] The opposite Schubert variety $X^v$ is smooth at the flag $wB$ if and only if
\[\sigma_v(w) = \prod_{\tiny \begin{array}{c} \alpha \in \Phi^+ s.t. \\ v \not \leq s_{\alpha}w \end{array}} \alpha.\]
\end{theorem}

\subsection{Restriction to fixed points and GKM theory} \label{section: GKM theory}
The two previous points---(1) constructing the cohomology of the flag variety and (2) the topology of Schubert varieties---fit together in a natural framework that also fall out  from Billey's formula: the family of polynomials associated to $\sigma_v$ by Billey's formula actually represent the equivariant Schubert class corresponding to $v$ and the individual polynomials $\sigma_v(w)$ encode topological information about the Schubert variety at the fixed point $w$.  This is now viewed as a part of a larger topological construction of torus-equivariant cohomology for a wider class of varieties that is often referred to as {\em GKM theory}.  (Andersen-Jantzen-Soergel's endomorphism algebras encapsulate key properties of GKM theory from a purely algebraic point of view.)

We describe GKM theory for flag varieties.  To begin, the inclusion $(G/B)^T \hookrightarrow G/B$ induces a map
\[H^*_T(G/B) \rightarrow H^*_T((G/B)^T).\]
For the flag variety (and many other varieties, including $G/P$), this map is an injection:
\[H^*_T(G/B) \hookrightarrow H^*_T((G/B)^T).\]

The equivariant cohomology of a point is
\[H^*_T(pt) \cong \mathbb{C}[\alpha_1, \alpha_2, \ldots, \alpha_n].\]
The fixed points $(G/B)^T$ are the flags associated to Weyl group elements $\{wB: w \in W\}$ so
\[H^*_T((G/B)^T) \cong \bigoplus_{w \in W}   \mathbb{C}[\alpha_1, \alpha_2, \ldots, \alpha_n].\]

Inclusion of fixed points induces an injection
\[H^*_T(G/B) \hookrightarrow \bigoplus_{w \in W}   \mathbb{C}[\alpha_1, \alpha_2, \ldots, \alpha_n].\]

This brings us to the key point about Billey's formula.

\begin{fact}The image of the Schubert class $[X^v]$ under this map is
\[ [X^v] \mapsto (\sigma_v(w))_{w \in W}\]
where $\sigma_v(w)$ are the very same polynomials given by Billey's formula.
\end{fact}

In fact, Kumar showed that Billey's formula holds for affine flag manifolds \cite[Appendix]{Bil99} and for any $G/P$ \cite{Kum02}.

When the flag $wB$ is a smooth point in the variety $X^v$ then the polynomials $\sigma_v(w)$ describe the (torus weight on the) tangent space to the variety $X^v$ at the fixed point $w$. Geometrically, then, Kumar's criterion says that a fixed point is smooth in a Schubert variety if Billey's formula is the product of certain distinct roots.

\begin{figure}[h]
\setlength{\unitlength}{0.01in}

\begin{picture}(150,130)(200,-45)
\thicklines
\put(270,-20){\color{red} \line(0,1){80}}
\put(270,-20){\color{green} \line(-1,1){20}}
\put(270,-20){\line(1,1){20}}

\put(250,0){\color{red} \line(0,1){40}}
\put(250,0){\line(1,1){40}}
\put(290,0){\color{red} \line(0,1){40}}
\put(290,0){\color{green} \line(-1,1){40}}

\put(250,40){\line(1,1){20}}
\put(290,40){\color{green} \line(-1,1){20}}

\put(270,60){\circle*{5}}
\put(250,40){\circle*{5}}
\put(290,40){\circle*{5}}
\put(290,0){\circle*{5}}
\put(250,0){\circle*{5}}
\put(270,-20){\circle*{5}}

\put(280,-25){$e$}
\put(226,0){$s_1$}
\put(212,40){$s_1s_2$}
\put(300,0){$s_2$}
\put(300,35){$s_2s_1$}
\put(255,70){$s_1s_2s_1$}
\put(185,-40){The Bruhat graph for $S_3$}
\end{picture}
\hspace{0.75in}
\begin{picture}(150,130)(210,-45)
\thicklines
\put(270,-20){\color{red} \line(0,1){80}}
\put(270,-20){\color{green} \line(-1,1){20}}
\put(270,-20){\line(1,1){20}}

\put(250,0){\color{red} \line(0,1){40}}
\put(250,0){\line(1,1){40}}
\put(290,0){\color{red} \line(0,1){40}}
\put(290,0){\color{green} \line(-1,1){40}}

\put(250,40){\line(1,1){20}}
\put(290,40){\color{green} \line(-1,1){20}}

\put(270,60){\circle*{5}}
\put(250,40){\circle*{5}}
\put(290,40){\circle*{5}}
\put(290,0){\circle*{5}}
\put(250,0){\circle*{5}}
\put(270,-20){\circle*{5}}

\put(280,-25){$0$}
\put(230,0){$\alpha_1$}
\put(230,40){$\alpha_1$}
\put(300,0){$0$}
\put(300,35){$\alpha_1+\alpha_2$}
\put(250,70){$\alpha_1+\alpha_2$}
\put(200,-40){The Schubert class $\sigma_{s_1}$}
\end{picture}
\caption{The Bruhat graph for $GL_3/B$ and a class in $H^*_T(GL_3/B)$}
\label{figure: example of schubert class}
\end{figure}
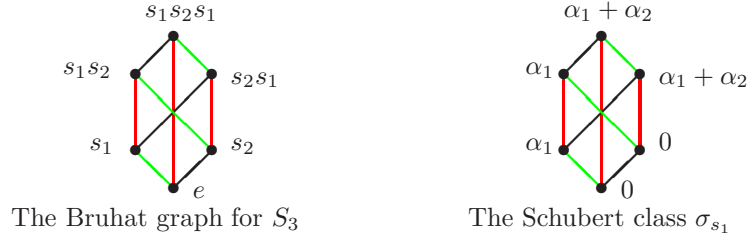

We often represent equivariant Schubert classes combinatorially.  The {\em Bruhat graph} of $W$ is the graph with vertices given by $w \in W$ and with edges $s_{\alpha}w \leftrightarrow w$ between vertices $w, s_{\alpha}w$ that differ by left-multiplication by a reflection.  We then describe the class $\sigma_v$  by  labeling each vertex $w$ of the Bruhat graph with the polynomial $\sigma_v(w)$ as shown in Figure \ref{figure: example of schubert class}.

\begin{remark}[GKM (Goresky-Kottwitz-MacPherson) theory] This construction of Schubert classes (and the cohomology ring) holds much more generally than just for the flag variety.  In fact, we can identify the image  $H^*_T(X) \hookrightarrow H^*_T(X^T)$ for a large family of varieties $X$.  In many cases---now referred to as GKM theory---the image can be described as a subset of polynomials labeling the vertices of a graph associated to $X$.  The conditions defining which polynomials are allowed are given by a straightforward combinatorial algorithm. GKM theory builds on work of Chang-Skjelbred, Kirwan, Atiyah-Bott, Guillemin-Sternberg, and many others \cite{ChaSkj74, Kir84, AtiBot84, GuiSte99}.
\end{remark}

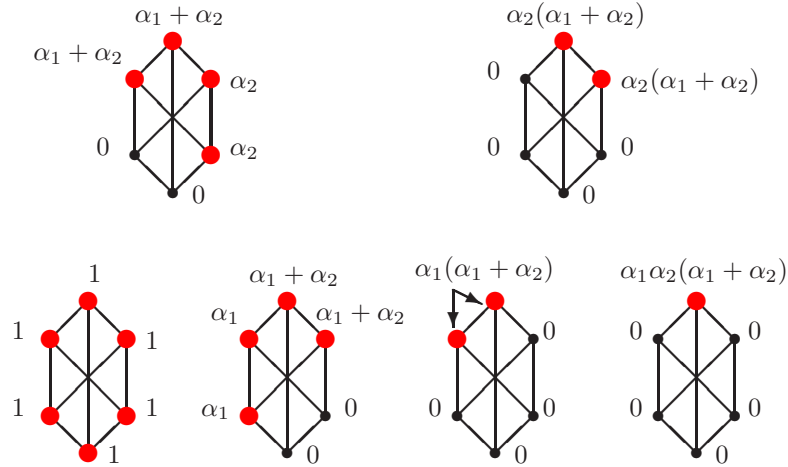
\begin{figure}[h]
\setlength{\unitlength}{0.01in}

\begin{picture}(150,125)(60,-40)
\thicklines

\put(100,-20){\line(0,1){80}}
\put(100,-20){\line(-1,1){20}}
\put(100,-20){\line(1,1){20}}

\put(80,0){\line(0,1){40}}
\put(80,0){\line(1,1){40}}
\put(120,0){\line(0,1){40}}
\put(120,0){\line(-1,1){40}}

\put(80,40){\line(1,1){20}}
\put(120,40){\line(-1,1){20}}

\put(100,60){\color{red} \circle*{9}}
\put(80,40){\color{red} \circle*{9}}
\put(120,40){\color{red} \circle*{9}}
\put(120,0){\color{red} \circle*{9}}
\put(80,0){\circle*{5}}
\put(100,-20){\circle*{5}}

\put(110,-25){$0$}
\put(60,0){$0$}
\put(27,50){$\alpha_1+\alpha_2$}
\put(130,0){$\alpha_2$}
\put(130,35){$\alpha_2$}
\put(80,70){$\alpha_1+\alpha_2$}
\end{picture} 
\begin{picture}(150,125)(130,-40)
\thicklines

\put(220,-20){\line(0,1){80}}
\put(220,-20){\line(-1,1){20}}
\put(220,-20){\line(1,1){20}}

\put(200,0){\line(0,1){40}}
\put(200,0){\line(1,1){40}}
\put(240,0){\line(0,1){40}}
\put(240,0){\line(-1,1){40}}

\put(200,40){\line(1,1){20}}
\put(240,40){\line(-1,1){20}}

\put(220,60){\color{red} \circle*{9}}
\put(200,40){\circle*{5}}
\put(240,40){\color{red} \circle*{9}}
\put(240,0){\circle*{5}}
\put(200,0){\circle*{5}}
\put(220,-20){\circle*{5}}

\put(230,-25){$0$}
\put(180,0){$0$}
\put(180,40){$0$}
\put(250,0){$0$}
\put(250,35){$\alpha_2(\alpha_1+\alpha_2)$}
\put(190,70){$\alpha_2(\alpha_1+\alpha_2)$}
\end{picture}

%\vspace{0.1in}

%\hspace{-0.25in}
\begin{picture}(80,125)(300,-30)
\thicklines

\put(320,-20){\line(0,1){80}}
\put(320,-20){\line(-1,1){20}}
\put(320,-20){\line(1,1){20}}

\put(300,0){\line(0,1){40}}
\put(300,0){\line(1,1){40}}
\put(340,0){\line(0,1){40}}
\put(340,0){\line(-1,1){40}}

\put(300,40){\line(1,1){20}}
\put(340,40){\line(-1,1){20}}

\put(320,60){\color{red} \circle*{9}}
\put(300,40){\color{red} \circle*{9}}
\put(340,40){\color{red} \circle*{9}}
\put(340,0){\color{red} \circle*{9}}
\put(300,0){\color{red} \circle*{9}}
\put(320,-20){\color{red} \circle*{9}}

\put(330,-25){$1$}
\put(280,0){$1$}
\put(280,40){$1$}
\put(350,0){$1$}
\put(350,35){$1$}
\put(320,70){$1$}
\end{picture}
\hspace{0.25in}
\begin{picture}(80,125)(110,-30)
\thicklines

\put(120,-20){\line(0,1){80}}
\put(120,-20){\line(-1,1){20}}
\put(120,-20){\line(1,1){20}}

\put(100,0){\line(0,1){40}}
\put(100,0){\line(1,1){40}}
\put(140,0){\line(0,1){40}}
\put(140,0){\line(-1,1){40}}

\put(100,40){\line(1,1){20}}
\put(140,40){\line(-1,1){20}}

\put(120,60){\color{red} \circle*{9}}
\put(100,40){\color{red} \circle*{9}}
\put(140,40){\color{red} \circle*{9}}
\put(140,0){\circle*{5}}
\put(100,0){\color{red} \circle*{9}}
\put(120,-20){\circle*{5}}

\put(130,-25){$0$}
\put(75,0){$\alpha_1$}
\put(80,50){$\alpha_1$}
\put(150,0){$0$}
\put(135,50){$\alpha_1+\alpha_2$}
\put(100,70){$\alpha_1+\alpha_2$}
\end{picture}
\hspace{0.25in}
\begin{picture}(60,125)(210,-30)
\thicklines

\put(320,-20){\line(0,1){80}}
\put(320,-20){\line(-1,1){20}}
\put(320,-20){\line(1,1){20}}

\put(300,0){\line(0,1){40}}
\put(300,0){\line(1,1){40}}
\put(340,0){\line(0,1){40}}
\put(340,0){\line(-1,1){40}}

\put(300,40){\line(1,1){20}}
\put(340,40){\line(-1,1){20}}

\put(320,60){\color{red} \circle*{9}}
\put(300,40){\circle*{5}}
\put(340,40){\circle*{5}}
\put(340,0){\circle*{5}}
\put(300,0){\circle*{5}}
\put(320,-20){\circle*{5}}

\put(330,-25){$0$}
\put(285,0){$0$}
\put(347,40){$0$}
\put(347,0){$0$}
\put(285,40){$0$}
\put(280,72){$\alpha_1\alpha_2(\alpha_1+\alpha_2)$}
\end{picture}
%\hspace{0.25in}
\begin{picture}(60,125)(380,-30)
\thicklines

\put(320,-20){\line(0,1){80}}
\put(320,-20){\line(-1,1){20}}
\put(320,-20){\line(1,1){20}}

\put(300,0){\line(0,1){40}}
\put(300,0){\line(1,1){40}}
\put(340,0){\line(0,1){40}}
\put(340,0){\line(-1,1){40}}

\put(300,40){\line(1,1){20}}
\put(340,40){\line(-1,1){20}}

\put(320,60){\color{red} \circle*{9}}
\put(300,40){\color{red} \circle*{9}}
\put(340,40){\circle*{5}}
\put(340,0){\circle*{5}}
\put(300,0){\circle*{5}}
\put(320,-20){\circle*{5}}

\put(330,-25){$0$}
\put(285,0){$0$}
\put(345,40){$0$}
\put(345,0){$0$}
%\put(260,50){\tiny $(t_1-t_2)(t_1-t_3)$}
\put(278,73){$\alpha_1(\alpha_1+\alpha_2)$}
\put(298, 66){\vector(2,-1){15}}
\put(298,66){\vector(0,-1){20}}
\end{picture}
\caption{Schubert classes for $H^*_T(GL_3/B)$}\label{figure: all Schubert classes}
\end{figure}

Figure \ref{figure: all Schubert classes} shows all of the Schubert classes for $GL_3/B$.

The central goal of Schubert calculus is to identify the products $\sigma_u \sigma_v$ in terms of Schubert classes, namely find explicit formulas for the coefficients $c_{uv}^w$ in the expansion
\[\sigma_u \sigma_v = \sum c_{uv}^w \sigma_w.\]
With the GKM presentation of the equivariant cohomology ring, several aspects of Schubert calculus become simpler:
\begin{itemize}
\item The ring structure of $H^*_T(G/B)$ is streamlined: we multiply and add these classes vertex-wise, using Billey's formula to identify the polynomials associated to each vertex.
\item The Schubert classes are ``upper-triangular" with respect to the partial order on $W$ in which $u \leq v$ if $u$ can be written as a subword of $v$.  (Figure \ref{figure: all Schubert classes} demonstrates this visually and Proposition \ref{proposition: upper-triangular} states it algebraically.)
\end{itemize}

For instance, Knutsen and Tao used GKM theory in their paper on puzzles to give a positive formula for equivariant structure constants in $H^*_T(G(k,n))$ \cite{KnuTao03}.  Their proof had three main steps: 1) proving that a small subset of structure constants determine all of the structure constants; 2) establishing certain instances of Billey's formula to identify this subset of structure constants; and 3) proving that the structure constants agreed in those instances with   the polynomials given by their puzzles. 

\begin{remark} Billey's formula holds in large part because of the strong combinatorial structures inherent in the flag variety.  Nonetheless, there has been some work to generalize this to a larger family of symplectic manifolds, notably by Goldin and Tolman \cite{GolTol09}.  Goldin and Tolman both generalize the concept of Schubert classes to a large family of symplectic manifolds and give a (generally-non-positive) formula for the restriction of these classes to the fixed points.  
\end{remark}

\section{Billey's formula for the Grassmannian $G(k,n)$}\label{section: applications}

In brief, Billey's formula turns calculations  in geometric Schubert calculus into combinatorics.  It can  be restricted to subvarieties of $G/P$, so Billey's formula can be used in a wide range of applications of Schubert calculus.  This section describes  applications for Grassmannians.  The next section shows how to use Billey's formula when GKM theory does not hold.

In particular cases, the combinatorics of Billey's formula can be made more explicit.  For instance, if the ambient variety is the Grassmannian $G(k,n)$ of $k$-dimensional subspaces of $\mathbb{C}^n$ then the combinatorial construction for Billey's formula is called {\em excited Young diagrams} by Ikeda and Naruse \cite{IkeNar09}, who discovered it independently after Kreiman \cite{Kre05}, or Knutson-Miller-Yong in the case of ordinary cohomology \cite{KMY09}.   

\begin{theorem}[Kreiman, Ikeda-Naruse] \label{theorem: excited young diagrams} In $G(k,n)$ the polynomial $\sigma_{\lambda}(\mu)$ is the sum of the excited Young diagrams of $\lambda$ inside $\mu$, weighted according to the position of the boxes.
\end{theorem}

An excited Young diagram is a Young diagram in which certain boxes are marked.  We draw our Young diagrams so that the partition $\mu = (\mu_1 \geq \mu_2 \geq \ldots)$ has $\mu_1$ boxes in the first row, $\mu_2$ boxes in the second row, and so on, with the rows aligned on the leftmost column.

\begin{defn}
A marked box in a Young diagram can be {\em excited} if it has empty boxes to the east, south, and southeast. When a box is excited, its marking moves to the box directly southeast. 
\end{defn}

We now assign a weight $(i,j)$ to each box in the Young diagram:
\begin{itemize}
\item $i$ is the index of the column in which the box is located (read from the left).
\item $j$ is (the number of rows at or below the box) $+$ (the index of the rightmost column in the box's row).
\end{itemize}

To construct all excited Young diagrams of $\lambda$ inside $\mu$, follow the resulting simple algorithm (which we demonstrate with an example):

\begin{enumerate}
\item Mark the boxes of $\lambda$ inside $\mu$.

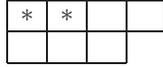
\[\begin{array}{|c|c|c|c|}
\cline{1-4} * & * &\hspace{0.5em} &\hspace{0.5em} \\
\cline{1-4} & & &\multicolumn{1}{|c}{} \\
\cline{1-3} \multicolumn{4}{c}{} \\
\end{array}\]

\item Excite all marked boxes in all possible ways.  

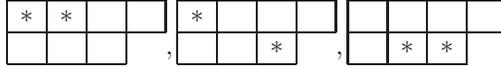
\[\begin{array}{|c|c|c|c|}
\cline{1-4} * & * & \hspace{0.5em} & \hspace{0.5em}\\
\cline{1-4} & & &\multicolumn{1}{|c}{} \\
\cline{1-3} \multicolumn{4}{c}{} \\
\end{array}, \begin{array}{|c|c|c|c|}
\cline{1-4} * & \hspace{0.5em} & \hspace{0.5em}& \hspace{0.5em} \\
\cline{1-4} & & * &\multicolumn{1}{|c}{} \\
\cline{1-3} \multicolumn{4}{c}{} \\
\end{array}
, \begin{array}{|c|c|c|c|}
\cline{1-4} \hspace{0.5em} & & & \hspace{0.5em}\\
\cline{1-4} & * & * &\multicolumn{1}{|c}{} \\
\cline{1-3} \multicolumn{4}{c}{} \\
\end{array}
\]

\item Sum the diagrams, weighted by position of marked boxes.

\[\begin{array}{c}
\tiny{\begin{array}{|c|c|c|c|}
\cline{1-4} 1,6 & 2,6 & \hspace{1em} & \hspace{1em}\\
\cline{1-4} & & &\multicolumn{1}{|c}{} \\
\cline{1-3} \multicolumn{4}{c}{} \\
\end{array}  + 
\begin{array}{|c|c|c|c|}
\cline{1-4} 1,6 & \hspace{1em} & \hspace{1em}& \hspace{1em} \\
\cline{1-4} & & 3,4 &\multicolumn{1}{|c}{} \\
\cline{1-3} \multicolumn{4}{c}{} \\
\end{array}
 + 
 \begin{array}{|c|c|c|c|}
\cline{1-4} \hspace{1em} & & & \hspace{1em}\\
\cline{1-4} & 2,4 & 3,4 &\multicolumn{1}{|c}{} \\
\cline{1-3} \multicolumn{4}{c}{} \\
\end{array}} = \\
(t_1-t_6)(t_2-t_6) + (t_1-t_6)(t_3-t_4) + (t_2-t_4)(t_3-t_4)\end{array}
\]
\end{enumerate}

We sketch the main ideas of the proof of Theorem \ref{theorem: excited young diagrams}, though we omit details.  Young diagrams represent fixed points in $G(k,n)$ as well as permutations.  These Grassmannian permutations satisfy strong algebraic conditions: they are the product of  descending strings of simple reflections $s_j s_{j-1} s_{j-2} \ldots s_{j-j'}$ and the rightmost reflection is $s_k$.  Excited Young diagrams enumerate all possible reduced subwords for $\lambda$ inside the word for $\mu$ using the combinatorics of $S_n$.

\section{Billey's formula for subvarieties} \label{section: poset pinball}

GKM theory cannot be applied to all varieties, not even if those varieties have rich geometric and combinatorial structures. We can nonetheless construct a kind of GKM theory for suitable subvarieties $Y$ of $G/B$ using Billey's formula.  Our strategy is to bootstrap information about Schubert classes to obtain a module basis for the equivariant cohomology of $Y$.

The subvariety $Y$ must satisfy two important conditions.
\begin{enumerate}
\item  $Y$ must admit the action of a (one-dimensional) subtorus $S \subseteq T$ under which the fixed points $Y^S \subseteq (G/B)^T$.
\item  $Y$ must be {\em equivariantly formal}.
\end{enumerate}
Subvarieties of $G/B$ that are defined by linear conditions, like Schubert varieties or Hessenberg varieties (see below) often satisfy the first condition.   Any variety with no odd-dimensional ordinary cohomology satisfies the second condition.

Here are the three tools that we have developed thus far:
\begin{itemize}
\item a graph $\Gamma$ built from $G/B$
\item fixed points $Y^S$ that are a subset of vertices in $\Gamma$
\item basis classes $\sigma_v$ of $H^*_T(G/B)$ that are indexed by the set of $v \in (G/B)^T$
\end{itemize}
The graph $\Gamma$ is the Bruhat graph (see Section \ref{section: GKM theory}, or Figure \ref{figure: example of schubert class} for an example), which is essential for GKM theory in a way that we did not need to make precise for our exposition.  The basis classes $\sigma_v$ are the Schubert classes, and they are specified by Billey's formula.  Finally, the fixed points $Y^S$ are guaranteed to be a subset of the vertices in $\Gamma$ by our hypotheses on $Y$.

We will use this information to identify a subset $\sigma_v$ whose images $\widetilde{\sigma}_v|_{Y^S}$ generate $H^*_{S}(Y)$.

Intuitively, our strategy is as follows:

\smallskip

\noindent $\bullet$ Drop one ball from each vertex $v \in Y^S$ in turn and let it roll down the edges of the graph $\Gamma$ until it lands at a vertex $r(v) \in (G/B)^T$
\[v \leadsto r(v)\]

This heuristic was originally called {\em poset pinball}, though others point out that {\em poset pachinko} may be more appropriate.  The fixed point $r(v)$ is called the {\em roll-down} vertex of $v$.

We first give the main theorem and then give several examples.  Essentially, the theorem says that the polynomials given by Billey's formula for the fixed points $v \in Y^S$ and the Schubert classes $\sigma_{r(v)}|_{Y^S}$ generate the $S$-equivariant cohomology ring of $Y$.

\begin{theorem}[Harada-Tymoczko] \label{theorem: roll-downs}
For each $u \in (G/B)^T$ let $\widetilde{\sigma}_u$ denote the image of the Schubert class $\sigma_u$ under the natural restriction map $H^*_T(G/B) \rightarrow \hspace{-0.7em} \rightarrow H^*_S(G/B)$. If $Y$ is an {\em appropriate} subspace of $G/B$ then 
\[ \left\{ \widetilde{\sigma}_{r(w)}(w') : w,w' \in Y^S \right\} \]
generates $H^*_S(Y)$ as a ring and forms a module basis for $H^*_S(Y)$.
\end{theorem}

\begin{figure}
\setlength{\unitlength}{0.01in}
\begin{picture}(150,125)(85,-40)
\thicklines

\put(100,-20){\line(0,1){80}}
\put(100,-20){\line(-1,1){20}}
\put(100,-20){\line(1,1){20}}

\put(80,0){\line(0,1){40}}
\put(80,0){\line(1,1){40}}
\put(120,0){\line(0,1){40}}
\put(120,0){\line(-1,1){40}}

\put(80,40){\line(1,1){20}}
\put(120,40){\line(-1,1){20}}

\put(100,60){\color{red} \circle*{7}}
\put(80,40){\color{red} \circle*{7}}
\put(120,40){\color{red} \circle*{7}}
\put(120,0){\color{red} \circle*{7}}
\put(80,0){\circle*{5}}
\put(100,-20){\circle*{5}}

%\put(110,-25){$0$}
%\put(60,0){$0$}
%\put(20,50){$t_1-t_3$}
\put(130,0){$t$}  %{$\alpha_2$}
\put(130,35){$t$}  %{$\alpha_2$}
\put(100,70){$t+t$}   %$\alpha_1+\alpha_2$}

\put(100,60){\color{red} \circle{12}}
\put(120,0){\color{red} \circle{12}}
\put(120,40){\color{red} \circle{12}}
\end{picture}
\begin{picture}(80,125)(315,-40)
\thicklines

\put(320,-20){\line(0,1){80}}
\put(320,-20){\line(-1,1){20}}
\put(320,-20){\line(1,1){20}}

\put(300,0){\line(0,1){40}}
\put(300,0){\line(1,1){40}}
\put(340,0){\line(0,1){40}}
\put(340,0){\line(-1,1){40}}

\put(300,40){\line(1,1){20}}
\put(340,40){\line(-1,1){20}}

\put(320,60){\color{red} \circle*{7}}
\put(300,40){\color{red} \circle*{7}}
\put(340,40){\color{red} \circle*{7}}
\put(340,0){\color{red} \circle*{7}}
\put(300,0){\color{red} \circle*{7}}
\put(320,-20){\color{red} \circle*{7}}

%\put(325,-25){$1$}
%\put(290,0){$1$}
%\put(280,40){$1$}
\put(350,0){$1$}
\put(350,35){$1$}
\put(320,70){$1$}

\put(320,60){\color{red} \circle{12}}
\put(340,40){\color{red} \circle{12}}
\put(340,0){\color{red} \circle{12}}

\end{picture}
\begin{picture}(80,125)(80,-40)
\thicklines

\put(120,-20){\line(0,1){80}}
\put(120,-20){\line(-1,1){20}}
\put(120,-20){\line(1,1){20}}

\put(100,0){\line(0,1){40}}
\put(100,0){\line(1,1){40}}
\put(140,0){\line(0,1){40}}
\put(140,0){\line(-1,1){40}}

\put(100,40){\line(1,1){20}}
\put(140,40){\line(-1,1){20}}

\put(120,60){\color{red} \circle*{7}}
\put(100,40){\color{red} \circle*{7}}
\put(140,40){\color{red} \circle*{7}}
\put(140,0){\circle*{5}}
\put(100,0){\color{red} \circle*{7}}
\put(120,-20){\circle*{5}}

%\put(130,-25){$0$}
%\put(70,-10){$t_1 - t_2$}
%\put(80,50){\tiny $t_1-t_2$}
\put(150,0){$0$}
\put(150,40){$t+t$}   %{$\alpha_1+\alpha_2$}
\put(125,70){$t+t$}   %{$\alpha_1+\alpha_2$}

\put(120,60){\color{red} \circle{12}}
\put(140,40){\color{red} \circle{12}}
\put(140,0){\color{red} \circle{12}}
\end{picture}

\caption{Restricted Schubert classes for the Springer variety}\label{figure: restricted Schubert classes}
\end{figure}
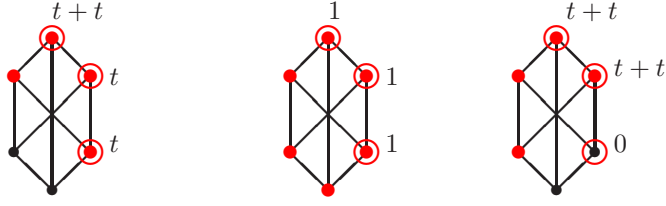

Figure \ref{figure: restricted Schubert classes} shows this in practice.  Each roll-down permutation gives a Schubert class; the vertices labeled by nonzero polynomials are emphasized, and the lowest of the marked vertices is the roll-down permutation.  The circled vertices are the elements of $Y^S$.  We record  Billey's formula at those vertices and discard all other information.  The resulting classes generate $H^*_S(Y)$ as a subring of $\left(\mathbb{C}[t]\right)^3$.

Various different sets of conditions for  $Y$ can be used in this theorem \cite{HarTym11, HarTym}.  Our main examples are subvarieties of $G/B$ though the results hold for a wider class of  ambient varieties than even $G/P$. 

Poset pinball applies to various important subvarieties of $G/B$, including both smooth and singular examples.

\smallskip

\noindent (1) {\bf Schubert varieties, $P/B$} \\
\noindent Poset pinball can be used for any subvariety to which regular GKM theory applies---like Schubert varieties and $P/B$. 

\smallskip

\noindent (2) {\bf Springer fibers of a nilpotent matrix} $X: \mathbb{C}^n \rightarrow \mathbb{C}^n$ 

\noindent The Springer variety of $X$ is the subvariety
\[Spr(X) = \left\{ \textup{flags }gB \in GL_n(\mathbb{C})/B : g^{-1}Xg \textup{ is upper-triangular}\right\}.\]
It is important in geometric representation theory: Springer showed that the symmetric group $S_n$ acts on the cohomology of the Springer variety \cite{Spr76}.   Moreover, the top-dimensional cohomology of $Spr(X)$ is the irreducible representation corresponding to the partition of $n$ given by the sizes of the Jordan blocks of $X$.  The combinatorics of Young tableaux give important information about the components and Betti numbers of Springer varieties \cite{Shi85, Spa76,  Tym06}.

\smallskip

\noindent (3) {\bf Peterson varieties for a regular nilpotent} $X \in \mathfrak{g}$

\noindent The Peterson variety of a regular nilpotent $X \in \mathfrak{g}$ is 
\[Pet = \left\{ \textup{flags }gB \in G/B : g^{-1}Xg \in \mathfrak{b} \oplus \bigoplus_{\alpha_i \in \Delta} \mathfrak{g}_{-\alpha_i}\right\}\]
Kostant showed that the quantum cohomology ring of the flag variety is isomorphic to the coordinate ring of a dense open subvariety of the Peterson variety \cite{Kos96}.  Rietsch gave an explicit isomorphism in which the quantum parameters correspond to certain determinants, giving a beautiful collection of determinantal identities \cite{Rie03}.  The geometry and topology of Peterson varieties is again deeply linked to combinatorics of permutations and partitions \cite{Pre, InsYon12}.

\smallskip

\noindent (4) {\bf Nilpotent Hessenberg varieties in type $A$ for a nilpotent $X \in \mathfrak{gl}_n$ and appropriate subspace $H \subseteq \mathfrak{gl}_n$}

\noindent Both Springer varieties and Peterson varieties are examples of a larger class of varieties called Hessenberg varieties.    Hessenberg varieties have two parameters: a nilpotent element $X \in \mathfrak{gl}_n$ and a subspace $H\subseteq \mathfrak{gl}_n$.  The subspace $H$ must satisfy two constraints: that $H$ contains the Borel subalgebra $\mathfrak{b}$ and that $H$ is closed under Lie bracket with $\mathfrak{b}$ in the sense that $[H, \mathfrak{b}] \subseteq H$.  The Hessenberg variety for these two parameters $X$ and $H$ is defined as
\[ Hess(X,H) = \left\{ \textup{flags }gB \in GL_n(\mathbb{C})/B : g^{-1}Xg \in H\right\}.\]
Hessenberg varieties have a  kind of cell decomposition enumerated by Young tableaux that relates the geometry  of the cells with the combinatorics of the tableaux \cite{Tym06}.

\subsection{Examples of poset pinball}  

We give two examples of the roll-down process. 

We use the subtorus $S \subseteq T$ that consists of diagonal matrices with $t^i$ in the $i^{th}$ row, for each $i$ with $1 \leq i \leq n$ and each $t \in \mathbb{C}$.  This means that the map $H^*_T(G/B) \rightarrow \hspace{-0.7em} \rightarrow H^*_S(G/B)$ sends each $\alpha_i \mapsto t$.

Both of our examples are for subregular Springer varieties in $GL_n(\mathbb{C})/B$.  In other words, the nilpotent matrix $X$ has two Jordan blocks of dimensions $n-1$ and $1$.

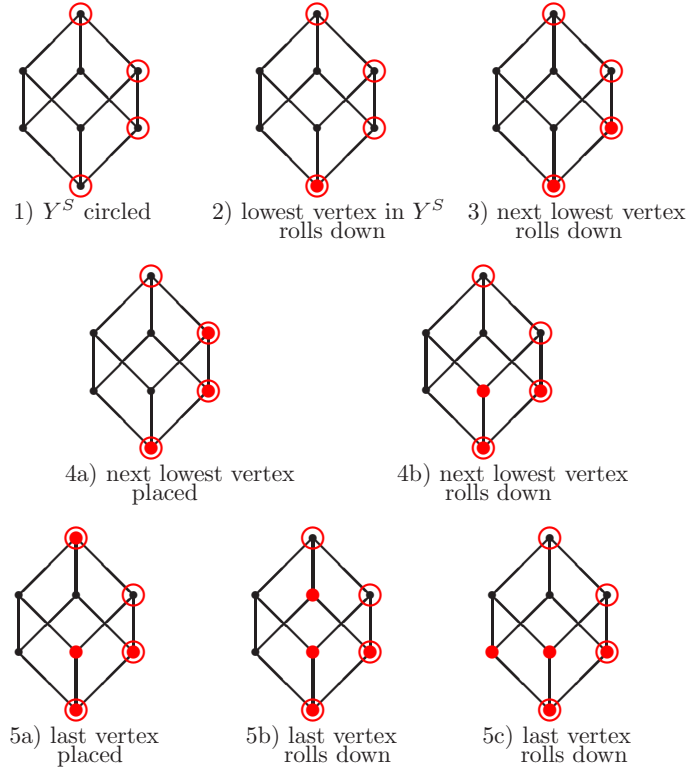
\begin{figure}[h]
\begin{center}
\setlength{\unitlength}{0.005in}
\begin{picture}(120,180)(-60,0)
\thicklines

\put(0,0){\circle*{9}}
\multiput(-60,60)(60,0){3}{\circle*{9}}
\multiput(-60,120)(60,0){3}{\circle*{9}}
\put(0,180){\circle*{9}}

\put(0,0){\line(-1,1){60}}
\put(0,0){\line(0,1){60}}
\put(0,0){\line(1,1){60}}

\multiput(-60,60)(120,0){2}{\line(0,1){60}}
\multiput(-60,60)(60,0){2}{\line(1,1){60}}
\multiput(0,60)(60,0){2}{\line(-1,1){60}}

\put(0,180){\line(-1,-1){60}}
\put(0,180){\line(0,-1){60}}
\put(0,180){\line(1,-1){60}}

\put(0,0){\color{red} \circle{22}}
\put(60,60){\color{red} \circle{22}}
\put(60,120){\color{red} \circle{22}}
\put(0,180){\color{red} \circle{22}}
\put(-70, -35){\small 1) $Y^S$ circled}
\end{picture}
\hspace{0.5in}
\setlength{\unitlength}{0.005in}
\begin{picture}(120,180)(-60,0)
\thicklines

\put(0,0){\circle*{9}}
\multiput(-60,60)(60,0){3}{\circle*{9}}
\multiput(-60,120)(60,0){3}{\circle*{9}}
\put(0,180){\circle*{9}}

\put(0,0){\line(-1,1){60}}
\put(0,0){\line(0,1){60}}
\put(0,0){\line(1,1){60}}

\multiput(-60,60)(120,0){2}{\line(0,1){60}}
\multiput(-60,60)(60,0){2}{\line(1,1){60}}
\multiput(0,60)(60,0){2}{\line(-1,1){60}}

\put(0,180){\line(-1,-1){60}}
\put(0,180){\line(0,-1){60}}
\put(0,180){\line(1,-1){60}}

\put(0,0){\color{red} \circle{22}}
\put(0,0){\color{red} \circle*{13}}
\put(60,60){\color{red} \circle{22}}
\put(60,120){\color{red} \circle{22}}
\put(0,180){\color{red} \circle{22}}
\put(-110, -35){\small 2) lowest vertex in $Y^S$}
\put(-40, -55){\small rolls down}
\end{picture}
\hspace{0.5in}
\setlength{\unitlength}{0.005in}
\begin{picture}(120,180)(-60,0)
\thicklines

\put(0,0){\circle*{9}}
\multiput(-60,60)(60,0){3}{\circle*{9}}
\multiput(-60,120)(60,0){3}{\circle*{9}}
\put(0,180){\circle*{9}}

\put(0,0){\line(-1,1){60}}
\put(0,0){\line(0,1){60}}
\put(0,0){\line(1,1){60}}

\multiput(-60,60)(120,0){2}{\line(0,1){60}}
\multiput(-60,60)(60,0){2}{\line(1,1){60}}
\multiput(0,60)(60,0){2}{\line(-1,1){60}}

\put(0,180){\line(-1,-1){60}}
\put(0,180){\line(0,-1){60}}
\put(0,180){\line(1,-1){60}}

\put(0,0){\color{red} \circle{22}}
\put(0,0){\color{red} \circle*{13}}
\put(60,60){\color{red} \circle{22}}
\put(60,60){\color{red} \circle*{13}}
\put(60,120){\color{red} \circle{22}}
\put(0,180){\color{red} \circle{22}}
\put(-90, -35){\small 3) next lowest vertex}
\put(-40, -55){\small rolls down}
\end{picture}
\end{center}

\vspace{0.4in}

\begin{center}
\setlength{\unitlength}{0.005in}
\begin{picture}(120,180)(-60,0)
\thicklines

\put(0,0){\circle*{9}}
\multiput(-60,60)(60,0){3}{\circle*{9}}
\multiput(-60,120)(60,0){3}{\circle*{9}}
\put(0,180){\circle*{9}}

\put(0,0){\line(-1,1){60}}
\put(0,0){\line(0,1){60}}
\put(0,0){\line(1,1){60}}

\multiput(-60,60)(120,0){2}{\line(0,1){60}}
\multiput(-60,60)(60,0){2}{\line(1,1){60}}
\multiput(0,60)(60,0){2}{\line(-1,1){60}}

\put(0,180){\line(-1,-1){60}}
\put(0,180){\line(0,-1){60}}
\put(0,180){\line(1,-1){60}}

\put(0,0){\color{red} \circle{22}}
\put(0,0){\color{red} \circle*{13}}
\put(60,60){\color{red} \circle{22}}
\put(60,60){\color{red} \circle*{13}}
\put(60,120){\color{red} \circle{22}}
\put(60,120){\color{red} \circle*{13}}
\put(0,180){\color{red} \circle{22}}
\put(-90, -35){\small 4a) next lowest vertex}
\put(-20, -55){\small placed}
\end{picture}
\hspace{1in}
\setlength{\unitlength}{0.005in}
\begin{picture}(120,180)(-60,0)
\thicklines

\put(0,0){\circle*{9}}
\multiput(-60,60)(60,0){3}{\circle*{9}}
\multiput(-60,120)(60,0){3}{\circle*{9}}
\put(0,180){\circle*{9}}

\put(0,0){\line(-1,1){60}}
\put(0,0){\line(0,1){60}}
\put(0,0){\line(1,1){60}}

\multiput(-60,60)(120,0){2}{\line(0,1){60}}
\multiput(-60,60)(60,0){2}{\line(1,1){60}}
\multiput(0,60)(60,0){2}{\line(-1,1){60}}

\put(0,180){\line(-1,-1){60}}
\put(0,180){\line(0,-1){60}}
\put(0,180){\line(1,-1){60}}

\put(0,0){\color{red} \circle{22}}
\put(0,0){\color{red} \circle*{13}}
\put(60,60){\color{red} \circle{22}}
\put(60,60){\color{red} \circle*{13}}
\put(60,120){\color{red} \circle{22}}
\put(0,60){\color{red} \circle*{13}}
\put(0,180){\color{red} \circle{22}}
\put(-90, -35){\small 4b) next lowest vertex}
\put(-40, -55){\small rolls down}
\end{picture}
\end{center}

\vspace{0.4in}
\begin{center}
\setlength{\unitlength}{0.005in}
\begin{picture}(120,180)(-60,0)
\thicklines

\put(0,0){\circle*{9}}
\multiput(-60,60)(60,0){3}{\circle*{9}}
\multiput(-60,120)(60,0){3}{\circle*{9}}
\put(0,180){\circle*{9}}

\put(0,0){\line(-1,1){60}}
\put(0,0){\line(0,1){60}}
\put(0,0){\line(1,1){60}}

\multiput(-60,60)(120,0){2}{\line(0,1){60}}
\multiput(-60,60)(60,0){2}{\line(1,1){60}}
\multiput(0,60)(60,0){2}{\line(-1,1){60}}

\put(0,180){\line(-1,-1){60}}
\put(0,180){\line(0,-1){60}}
\put(0,180){\line(1,-1){60}}

\put(0,0){\color{red} \circle{22}}
\put(0,0){\color{red} \circle*{13}}
\put(60,60){\color{red} \circle{22}}
\put(60,60){\color{red} \circle*{13}}
\put(60,120){\color{red} \circle{22}}
\put(0,60){\color{red} \circle*{13}}
\put(0,180){\color{red} \circle{22}}
\put(0,180){\color{red} \circle*{13}}
\put(-70, -35){\small 5a) last vertex}
\put(-20, -55){\small placed}
\end{picture}
\hspace{0.5in}
\setlength{\unitlength}{0.005in}
\begin{picture}(120,180)(-60,0)
\thicklines

\put(0,0){\circle*{9}}
\multiput(-60,60)(60,0){3}{\circle*{9}}
\multiput(-60,120)(60,0){3}{\circle*{9}}
\put(0,180){\circle*{9}}

\put(0,0){\line(-1,1){60}}
\put(0,0){\line(0,1){60}}
\put(0,0){\line(1,1){60}}

\multiput(-60,60)(120,0){2}{\line(0,1){60}}
\multiput(-60,60)(60,0){2}{\line(1,1){60}}
\multiput(0,60)(60,0){2}{\line(-1,1){60}}

\put(0,180){\line(-1,-1){60}}
\put(0,180){\line(0,-1){60}}
\put(0,180){\line(1,-1){60}}

\put(0,0){\color{red} \circle{22}}
\put(0,0){\color{red} \circle*{13}}
\put(60,60){\color{red} \circle{22}}
\put(60,60){\color{red} \circle*{13}}
\put(60,120){\color{red} \circle{22}}
\put(0,60){\color{red} \circle*{13}}
\put(0,180){\color{red} \circle{22}}
\put(0,120){\color{red} \circle*{13}}
\put(-70, -35){\small 5b) last vertex}
\put(-30, -55){\small rolls down}
\end{picture}
\hspace{0.5in}
\setlength{\unitlength}{0.005in}
\begin{picture}(120,180)(-60,0)
\thicklines

\put(0,0){\circle*{9}}
\multiput(-60,60)(60,0){3}{\circle*{9}}
\multiput(-60,120)(60,0){3}{\circle*{9}}
\put(0,180){\circle*{9}}

\put(0,0){\line(-1,1){60}}
\put(0,0){\line(0,1){60}}
\put(0,0){\line(1,1){60}}

\multiput(-60,60)(120,0){2}{\line(0,1){60}}
\multiput(-60,60)(60,0){2}{\line(1,1){60}}
\multiput(0,60)(60,0){2}{\line(-1,1){60}}

\put(0,180){\line(-1,-1){60}}
\put(0,180){\line(0,-1){60}}
\put(0,180){\line(1,-1){60}}

\put(0,0){\color{red} \circle{22}}
\put(0,0){\color{red} \circle*{13}}
\put(60,60){\color{red} \circle{22}}
\put(60,60){\color{red} \circle*{13}}
\put(60,120){\color{red} \circle{22}}
\put(0,60){\color{red} \circle*{13}}
\put(0,180){\color{red} \circle{22}}
\put(-60,60){\color{red} \circle*{13}}
\put(-70, -35){\small 5c) last vertex}
\put(-30, -55){\small rolls down}
\end{picture}
\vspace{0.2in}
\end{center}

\caption{Poset pinball on a Springer variety}\label{figure: springer variety example}
\end{figure}

Figure \ref{figure: springer variety example} shows the case when the subvariety $Y \subseteq GL_4(\mathbb{C})/B$ is the Springer variety associated to the matrix $X$ with $1$ in entries $(1,2), (2,3)$ and zero elsewhere.  (We only show relevant vertices and edges of the Bruhat graph $\Gamma$.)  The bold vertices in Figure \ref{figure: springer variety example} are the roll-down vertices; the circled vertices are the fixed points $Y^S$.

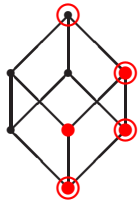
\begin{figure}[h]
\begin{center}
\setlength{\unitlength}{0.005in}
\begin{picture}(120,180)(-60,0)
\thicklines

\put(0,0){\circle*{9}}
\multiput(-60,60)(60,0){3}{\circle*{9}}
\multiput(-60,120)(60,0){3}{\circle*{9}}
\put(0,180){\circle*{9}}

\put(0,0){\line(-1,1){60}}
\put(0,0){\line(0,1){60}}
\put(0,0){\line(1,1){60}}

\multiput(-60,60)(120,0){2}{\line(0,1){60}}
\multiput(-60,60)(60,0){2}{\line(1,1){60}}
\multiput(0,60)(60,0){2}{\line(-1,1){60}}

\put(0,180){\line(-1,-1){60}}
\put(0,180){\line(0,-1){60}}
\put(0,180){\line(1,-1){60}}

\put(0,0){\color{red} \circle{22}}
\put(0,0){\color{red} \circle*{13}}
\put(60,60){\color{red} \circle{22}}
\put(60,60){\color{red} \circle*{13}}
\put(60,120){\color{red} \circle{22}}
\put(0,60){\color{red} \circle*{13}}
\put(0,180){\color{red} \circle{22}}
\put(60,120){\color{red} \circle*{13}}
\end{picture}
\end{center}

\caption{Alternate roll-downs}\label{figure: springer variety alternate example}
\end{figure}

This process is {\bf not} deterministic: for instance, instead of rolling down in step 5b of Figure \ref{figure: springer variety example}, the ball could instead have rolled right to give the alternate set of roll-down vertices shown in Figure \ref{figure: springer variety alternate example}.

\begin{figure}[h]
\begin{center}
\setlength{\unitlength}{0.01in}
{\begin{picture}(80,85)(-20,-20)

\thicklines

\put(17,-20){ \line(0,1){80}}
\put(-3,0){ \line(0,1){40}}
\put(37,0){ \line(0,1){40}}

\put(0,40){\line(1,1){20}}
\put(0,0){\line(1,1){40}}
\put(20,-20){\line(1,1){20}}

\put(38,42){\line(-1,1){20}}
\put(37,0){ \line(-1,1){40}}
\put(17,-20){ \line(-1,1){20}}

\put(20,60){\circle*{5}}
\put(0,40){\circle*{5}}
\put(40,40){\circle*{5}}
\put(0,0){\circle*{5}}
\put(40,0){\circle*{5}}
\put(20,-20){\circle*{5}}

\put(20,60){\color{red} \circle{12}}
\put(40,40){\color{red} \circle{12}}
\put(40,0){\color{red} \circle{12}}

\put(20,60){\color{red} \circle*{8}}
\put(40,40){\color{red} \circle*{8}}
\put(40,0){\color{red} \circle*{8}}

\put(47,-2){$s_2$}
\put(47,38){$s_2s_1$}
\put(27,58){$s_1s_2s_1$}

\put(-15,-35){\small Before rolling}
\end{picture}}
\hspace{0.25in}
\begin{picture}(60,45)(-5,-15)
\put(-10,45){\vector(1,0){110}}
\put(-10,30){Roll down from $s_2$,}
\put(-5,15){then $s_2s_1$,}
\put(-5,0){and then $s_1s_2s_1$}
\put(-10,-5){\vector(1,0){110}}
\end{picture}
\hspace{0.35in}
\setlength{\unitlength}{0.01in}
{\begin{picture}(80,85)(-20,-20)

\thicklines

\put(17,-20){ \line(0,1){80}}
\put(-3,0){ \line(0,1){40}}
\put(37,0){ \line(0,1){40}}

\put(0,40){\line(1,1){20}}
\put(0,0){\line(1,1){40}}
\put(20,-20){\line(1,1){20}}

\put(38,42){\line(-1,1){20}}
\put(37,0){ \line(-1,1){40}}
\put(17,-20){ \line(-1,1){20}}

\put(20,60){\circle*{5}}
\put(0,40){\circle*{5}}
\put(40,40){\circle*{5}}
\put(0,0){\circle*{5}}
\put(40,0){\circle*{5}}
\put(20,-20){\circle*{5}}

\put(20,60){\color{red} \circle{12}}
\put(40,40){\color{red} \circle{12}}
\put(40,0){\color{red} \circle{12}}

\put(20,-20){\color{red} \circle*{8}}
\put(0,0){\color{red} \circle*{8}}
\put(40,0){\color{red} \circle*{8}}
\put(-25,-35){\small After rolling down}

\end{picture}}
%\vspace{0.35in}
\end{center}
\caption{Placement of balls before and after rolling down}\label{figure: subregular springer for GL_3}
\end{figure}
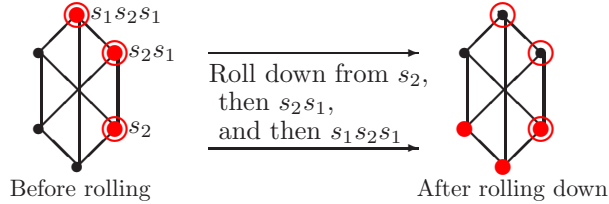

The next example is similar but in $GL_3(\mathbb{C})/B$.  The Springer variety is given by the matrix whose only nonzero entry is in position $(1,2)$.  We skip steps of the roll-down process because in this case it is essentially unique.  As long as the balls drop as far down as possible first from $s_2$, second from $s_2s_1$, and last from $s_1s_2s_1$, then the set of roll-downs must be the bottom three vertices, as  in Figure \ref{figure: subregular springer for GL_3}. 

%\begin{theorem} (Harada-T) Under appropriate circumstances
%\[ \langle \widetilde{\sigma}_{r(w)}(w') : w,w' \in Y^S\rangle \textup{ is } H^*_S(Y).\]
%\end{theorem}

Returning to Figure \ref{figure: restricted Schubert classes}, we see a pictorial summary of Theorem \ref{theorem: roll-downs}.  The three roll-downs from Figure \ref{figure: subregular springer for GL_3} are vertices $s_2$, $e$, and $s_1$.  Figure \ref{figure: all Schubert classes} shows the corresponding Schubert classes (as well as all other Schubert classes) for the flag variety; Billey's formula gives the polynomials at each vertex. Figure \ref{figure: subregular springer for GL_3} shows the classes for vertices $s_2$, $e$, and $s_1$ from left-to-right; now we only show polynomials at circled vertices, namely, the fixed points associated to the subvariety.

In short, the roll-down process selects certain Schubert classes for each subvariety of the flag variety.  Theorem \ref{theorem: roll-downs} then says that restricting those Schubert classes just to the fixed points in the subvariety---and discarding all additional information from Billey's formula---generates the equivariant cohomology of the subvariety.  Again, Billey's formula gives an explicit description of each basis class.

\subsection{Schubert calculus and poset pinball}\label{section: Schubert calculus for Petersons}  We use these results to actually compute Schubert calculus for suitable subvarieties of the flag variety, in this case for the Peterson variety $Y$ of type $A_{n-1}$.  First we identify roll-down classes for the Peterson variety.  Then we give two key formulas that together describe the entire cohomology ring of Peterson varieties: a Chevalley-Monk formula, which shows how to multiply any roll-down class by one of a specific family of generators; and a Giambelli formula, which expresses  an arbitrary roll-down class in terms of the generators.

For each subset $\mathcal{A} \subseteq \{1,2,\ldots,n-1\}$ define $v_{\mathcal{A}} = \prod_{i \in \mathcal{A}} s_i$.

\begin{theorem}\cite{HarTym11, HarTym}
Let $Y$ be the Peterson variety of type $A_{n-1}$ and $S$ be the subtorus of diagonal matrices with $t^i$ in the $i^{th}$ row for each $i$ with $1 \leq i \leq n$ and each $t \in \mathbb{C}$.

The equivariant cohomology $H^*_S(Y)$ is generated by
\[H^*_S(Y) = \left\langle \widetilde{\sigma}_{v_{\mathcal{A}}}(w):  \mathcal{A} \subseteq \{1,2,\ldots,n-1\} , wB \in Y^S \right\rangle.\]
\end{theorem}

This result is true in type $A_{n-1}$ \cite{HarTym11} and in general Lie type \cite{HarTym}.

We now give the Chevalley-Monk and Giambelli formulas for Peterson varieties.  The integers in these formulas are straightforward, but their precise descriptions  require distracting notation.

\begin{theorem} 
Let $Y$ be the Peterson variety of type $A_{n-1}$ and $S$ be the subtorus of diagonal matrices with $t^i$ in the $i^{th}$ row for each $i$ with $1 \leq i \leq n$ and each $t \in \mathbb{C}$.

\begin{enumerate}
\item {\bf Chevalley-Monk formula for Peterson varieties} \cite{HarTym11}
\[\widetilde{\sigma}_{s_{i}} \widetilde{\sigma}_{v_{\mathcal{A}}} = 
  t c_{i,\mathcal{A}}^\mathcal{A} \widetilde{\sigma}_{v_{\mathcal{A}}} + 
     \sum c_{i,\mathcal{A}}^\mathcal{B} \widetilde{\sigma}_{v_{\mathcal{B}}}\]
where the sum is over $\mathcal{B} \supseteq \mathcal{A}$ with $|\mathcal{B}| = |\mathcal{A}|+1$.  

\item {\bf Giambelli's formula for Peterson varieties} \cite{BayHar12}
\[\widetilde{\sigma}_{v_{\mathcal{A}}}  = c \prod_{i \in \mathcal{A}} \widetilde{\sigma}_{s_{i}}.\]
\end{enumerate}
All coefficients are explicit, easily-computed, positive integers.
\end{theorem}

Drellich recently extended this  to Peterson varieties of all Lie types \cite{Dre}.  She proves her results in several cases; however, the uniformity across type suggests an underlying topological cause. 

\section{Conjectures and open questions} \label{section: questions}

First we ask  for something like Billey's formula for varieties other than $G/P$.  We want a {\em positive} formula that {\em explicitly} describes each polynomial in the GKM presentation of $H^*_T(X)$.

\begin{question}
We ask for explicit, positive Billey formulas for varieties other than $G/P$.
\end{question}

Using the methods of Section \ref{section: poset pinball}, we can identify a Schubert variety $X_H$ for each regular nilpotent Hessenberg variety $Y_H$ so that

\[\{ \widetilde{\sigma}_{v}(w): vB, wB \in X_H \} \textup{ generates } H^*_{S}(X_H)\]
and we conjecture
\[\{ \widetilde{\sigma}_{v}(w): vB \in X_H, wB \in Y_H \} \textup{ generates } H^*_{S}(Y_H)\]
Analogous conjectures hold for Springer and all other nilpotent Hessenberg varieties.

\begin{question}Why?
\end{question}

In particular, is there a geometric reason for this similarity---e.g. a deformation of Hessenberg varieties to unions of Schubert varieties?

\end{document}